\documentclass[12pt]{article}
\usepackage{amsmath}
\usepackage{amssymb}
\usepackage{vatola}

\allowdisplaybreaks[4]
\def\q{\quad}

\def\qtq#1{\q\t{#1}\q}
\def\mod#1{\ (\text{\rm mod}\ #1)}
\def\t{\text}
\def\f{\frac}
\def\e{\equiv}
\def\b{\binom}
\def\v2{\vskip0.2cm}

\def\sls#1#2{(\f{#1}{#2})}
 \def\ls#1#2{\big(\f{#1}{#2}\big)}
\def\Ls#1#2{\Big(\f{#1}{#2}\Big)}

\let \pro=\proclaim
\let \endpro=\endproclaim

\begin{document}
 \centerline {\bf
Congruences for the Ap\'ery numbers modulo $p^3$}
\par\q\newline \centerline{Zhi-Hong Sun}
\newline \centerline{School of Mathematics
and Statistics}\centerline{Huaiyin Normal University}
\centerline{Huaian, Jiangsu 223300, P.R. China} \centerline{Email:
zhsun@hytc.edu.cn} \centerline{Homepage:
http://maths.hytc.edu.cn/szh1.htm}
\par\q\newline \centerline{}
\par\q
\par {\it Abstract. }
 Let $\{A'_n\}$ be the Ap\'ery numbers given by $A'_n=\sum_{k=0}^n\binom
 nk^2\binom{n+k}k.$ For any prime $p\equiv 3\pmod 4$ we show that $A'_{\frac{p-1}2}\equiv \frac{p^2}3\binom{\frac{p-3}2}{\frac{p-3}4}^{-2}\pmod {p^3}$. Let $\{t_n\}$ be given by
$$t_0=1,\ t_1=5\quad\hbox{and}\quad t_{n+1}=(8n^2+12n+5)t_n-4n^2(2n+1)^2t_{n-1}\ (n\ge 1).$$ We also obtain the congruences for $t_p\pmod {p^3},\ t_{p-1}\pmod {p^2}$ and $t_{\frac{p-1}2}\pmod {p^2}$, where $p$ is an odd prime.
 \par\q

\section*{1. Introduction}

\par \par For $s>1$ let $\zeta(s)=\sum_{n=1}^{\infty}\f 1{n^s}$.
In 1979, in order to prove that $\zeta(3)$ and
 $\zeta(2)$ are irrational,
 Ap\'ery [Ap] introduced the Ap\'ery numbers $\{A_n\}$ and $\{A'_n\}$
 given by
 $$ A_n= \sum_{k=0}^n\binom nk^2\binom{n+k}k^2\qtq{and}A'_n=\sum_{k=0}^n\b
 nk^2\b{n+k}k.$$
 The first few values of $A_n$ and $A'_n$ are shown below:
 $$\align &A_1=5,\ A_2=73,\ A_3=1445,\ A_4=33001,\ A_5=819005,\ A_6=21460825,
 \\&A'_1=3,\ A'_2=19,\ A'_3=147,\ A'_4=1251,\ A'_5=11253,\ A'_6=104959.
 \endalign$$
 It is well known (see [B2]) that
$$\align &(n+1)^3A_{n+1}=(2n+1)(17n(n+1)+5)A_n-n^3A_{n-1}\q (n\ge 1),
\\&(n+1)^2A'_{n+1}=(11n(n+1)+3)A'_n+n^2A'_{n-1}\q (n\ge 1).\endalign$$

\par Let $\Bbb Z^+$ denote the set of positive integers. In [B1] Beukers showed that for any prime $p>3$ and $m,r\in\Bbb Z^+$,
$$A_{mp^r-1}\e A_{mp^{r-1}-1}\mod {p^{3r}},\q
A'_{mp^r-1}\e A'_{mp^{r-1}-1}\mod {p^{3r}}.$$
Recently, Liu[L] showed that for any prime $p>3$ and $m,r\in\Bbb Z^+$,
$$\align &A_{mp^r}\e A_{mp^{r-1}}+\f 23C_mp^{3r}B_{p-3}
\mod {p^{3r+1}},
\\&A'_{mp^r}\e A'_{mp^{r-1}}+\f 13C'_mp^{3r}B_{p-3}\mod {p^{3r+1}},\endalign $$
where $\{B_n\}$ are the Bernoulli numbers given by $B_0=1$ and $\sum_{k=0}^{n-1}\b nkB_k=0\ (n\ge 2)$, 
$$C_m=\sum_{k=0}^m\b mk^2\b{m+k}k^2((m-k)^2-2km^2)\tag 1.1
$$
and
$$C'_m=\sum_{k=0}^m\b mk\b{m+k}k\big(2(m-k)^2-3m^2(m-k)-2k^2m\big).\tag 1.2$$
In [B2] Beukers conjectured that for any odd prime $p$,
$$A_{\f{p-1}2}'\e\cases 4x^2-2p\mod{p^2}&
\t{if $4\mid p-1$ and so $p=x^2+4y^2\ (x,y\in\Bbb Z)$,}
\\0\mod{p^2}&\t{if $p\e 3\mod 4$.}
\endcases\tag 1.3$$
This was proved by several authors including Ishikawa[I] ($p\e
1\pmod 4$), Van Hamme[VH]($p\e 3\pmod 4$) and Ahlgren[Ah].
\par In this paper, we prove that
$$A'_{\f{p-1}2}\e \f{p^2}3\b{\f{p-3}2}{\f{p-3}4}^{-2}\mod {p^3}\q\t{for any prime $p\e 3\mod 4$,}\tag 1.4$$
which was conjectured by the author in [Su5].
\par Let $\{t_n\}$ be given by
$$t_0=1,\ t_1=5\ \t{and}\ t_{n+1}=(8n^2+12n+5)t_n-4n^2(2n+1)^2t_{n-1}\ (n\ge 1).\tag 1.5$$
Then
$$t_1=5,\ t_2=89,\ t_3=3429,\ t_4=230481,\ t_5=23941125,\ t_6=3555578025.$$
In the paper, we investigate the identities and congruences for $\{t_n\}$. In particular, we show that for $n=0,1,2,\ldots$,  $$t_n^2=-(2n+1)!^2\sum_{k=0}^{2n+1}\b{2n+1+k}{2k}\b{2k}k^2\f 1{(-4)^k}
\sum_{i=1}^k\f 1{(2i-1)^2},$$
and obtain the congruences for $t_p\mod {p^3},\ t_{p-1}\mod {p^2},\ t_{\f{p\pm 1}2}\mod {p^2}$ and $t_{\f{p-3}4}\mod p$ for $p\e 3\mod 4$, where $p$ is an odd prime. For example, we have
$$t_p\e \big(1+4(-1)^{\f{p-1}2}\big)p^2\mod {p^3}\q\t{and}\q t_{\f{p-1}2}\e  pB_{p-1}-p+2^{p-1}-1\mod {p^2}.$$ 
\par Throughout this paper, the harmonic numbers $\{H_n\}$ are given by
$H_0=0$ and $H_n=1+\f 12+\cdots+\f 1n\ (n\ge 1)$, the Fermat quotient $q_p(a)=(a^{p-1}-1)/p$ and the Euler numbers $\{E_n\}$ are defined by
$$E_{2n-1}=0,\q E_0=1\q\t{and}\q E_{2n}=-\sum_{k=1}^{n} \b
{2n}{2k}E_{2n-2k}\q(n\ge 1).$$

\section*{2. Congruences for $A'_{p-1}$ modulo $p^3$}
\par\q Let $\{D_n\}$ be defined by $D_0=0$ and
$$D_n=2\sum_{1\le i<j\le n}\f 1{(2i-1)(2j-1)}=\Big(\sum_{i=1}^n\f 1{2i-1}\Big)^2-\sum_{i=1}^n\f 1{(2i-1)^2}\q(n\ge 1).$$
\pro{Lemma 2.1} For $n=0,1,2,\ldots$ we have
$$\sum_{k=0}^n\b nk(-1)^k\f{\b{2k}k}{4^k}D_k=\f{\b{2n}n}{4^n}D_n.$$
\endpro
\par{\it Proof.} Let
$$S_1(n)=\sum_{k=0}^n\b nk(-1)^k\f{\b{2k}k}{4^k}D_k\q\t{and}\q S_2(n)=\f{\b{2n}n}{4^n}D_n.$$
Then $S_0(0)=0=S_2(0),\ S_1(1)=0=S_2(1),\ S_1(2)=\f 14=S_2(2).$
Using the software Sigma we find that for $i=1,2$,
$$\align &8(n+1)(n+2)(n+3)S_i(n+3)-12(n+1)(n+2)(2n+3)S_i(n+2)
\\&+2(n+1)(12n^2+24n+13)S_i(n+1)
-(2n+1)^3S_i(n)=0\q(n=0,1,2,\ldots).\endalign$$
Hence, $S_1(n)=S_2(n)$ for $n=0,1,2,\ldots$. This proves the lemma.

\v2
\pro{Lemma 2.2 ([Su3, Theorem 2.2])} If $\{a_n\}$ is a sequence satisfying
$$\sum_{k=0}^n\b nk(-1)^ka_k=a_n\q(n=0,1,2,\ldots),$$
then
$$\sum_{k=0}^n\b nk\b{n+k}k(-1)^k a_k=0\q(n=1,3,5,\ldots).$$
\v2
\pro{Lemma 2.3} Let $p$ be an odd prime. Then
$$\align A'_{\f{p-1}2}&\e 1+\sum_{k=1}^{(p-1)/2}\f{\b{2k}k^3}{64^k}
\Big(1-p\sum_{i=1}^k\f 1{2i-1}\\&\q+\f{p^2}2\Big(\Big(\sum_{i=1}^k\f 1{2i-1}\Big)^2-3\sum_{i=1}^k\f 1{(2i-1)^2}\Big)\Big)\mod {p^3}.
\endalign$$
\endpro
\par{\it Proof.} Clearly
$$\align &A'_{\f{p-1}2}\\&=1+\sum_{k=1}^{\f{p-1}2}\f{\sls{p-1}2^2(\f{p-1}2-1)^2\cdots(\f{p-1}2-k+1)^2}{k!^2}
\cdot\f{(\f{p-1}2+1)(\f{p-1}2+2)\cdots(\f{p-1}2+k)}{k!}
\\&=1+\sum_{k=1}^{\f{p-1}2}
\f{(p-1)(p-3)\cdots(p-(2k-1))\cdot(p^2-1^2)(p^2-3^2)\cdots(p^2-(2k-1)^2)}{2^{3k}\cdot k!^3}
\\&\e 1+\sum_{k=1}^{(p-1)/2}\f{(1\cdot 3\cdot\cdots (2k-1))^3}{2^{3k}\cdot k!^3}
\Big(1-p\sum_{i=1}^k\f 1{2i-1}+\f{p^2}2D_k\Big)\\&\q\times\Big(1-p^2\sum_{i=1}^k\f 1{(2i-1)^2}\Big)
\\&\e 1+\sum_{k=1}^{(p-1)/2}\f{\b{2k}k^3}{64^k}
\Big(1-p\sum_{i=1}^k\f 1{2i-1}\\&\q+\f{p^2}2\Big(\Big(\sum_{i=1}^k\f 1{2i-1}\Big)^2-3\sum_{i=1}^k\f 1{(2i-1)^2}\Big)\Big)\mod {p^3}.
\endalign$$

\par\textbf{Lemma 2.4 ([Su6, Theorem 4.1])} {\sl Let $p$ be an odd
prime. Then
$$\sum_{k=0}^{p-1}\f{\b{2k}k^3}{64^k}
\e\begin{cases} 4x^2-2p-\f{p^2}{4x^2}\mod {p^3} &\t{if
$p=x^2+4y^2\e 1\mod 4$,}
\\-\f{p^2}4\b{(p-3)/2}{(p-3)/4}^{-2}\mod {p^3}
&\t{if $p\e 3\mod 4$}
\end{cases}$$}

\par  For an odd prime $p$ and rational $p-$integer $x$,
 the $p$-adic Gamma function $\Gamma_p(x)$ is
  defined by
  $$\Gamma_p(0)=1,\q \Gamma_p(n)=(-1)^n\prod_{
  \substack{k\in\{1,2,\ldots,n-1\}\\p\nmid
  k}}k\qtq{for}n=1,2,3,\ldots$$
  and
  $$\Gamma_p(x)=\lim_{\substack{n\in\{0,1,\ldots\}\\ |x-n|_p\rightarrow
  0}}
  \Gamma_p(n).$$

\pro{Lemma 2.5 ([T3, (9)])} Let $p$ be an odd prime. Then
$$\Gamma_p\Ls 14^4\e \begin{cases}
-\f 1{2^{p-1}}\b{\f{p-1}2}{\f{p-1}4}^2\Big(1-\f{p^2}2E_{p-3}\Big)
\pmod {p^3}&\t{if $4\mid p-1$,}
\\2^{p-3}(16+32p+(48-8E_{p-3})p^2)\b{\f{p-3}2}{\f{p-3}4}^{-2}
\pmod {p^3}&\t{if $4\mid p-3$.}
\end{cases}$$

\pro{Lemma 2.6 ([Su4, Theorem 2.8])} Let $p$ be a prime of the form $4k+1$ and so $p=x^2+4y^2$ with $x,y\in\Bbb Z$. Then
$$\f 1{2^{p-1}}\b{\f{p-1}2}{\f{p-1}4}^2\Big(1-\f{p^2}2E_{p-3}\Big)
\e 4x^2-2p-\f{p^2}{4x^2}\pmod {p^3}.$$

\pro{Lemma 2.7 ([T2])} For any prime $p>3$, 
$$\sum_{k=1}^{p-1}\f{\b{2k}k^3}{64^k}\sum_{i=1}^k\f 1{2i-1}
\e\cases 0\mod {p^2}&\t{if $p\e 1\mod 4$,}
\\-\f p{12}\Gamma_p\ls 14^4\mod {p^2}&\t{if $p\e 3\mod 4$}
\endcases$$
and
$$\sum_{k=1}^{p-1}\f{\b{2k}k^3}{64^k}\sum_{i=1}^k\f 1{(2i-1)^2}
\e\cases \f 12\Gamma_p\sls 14^4E_{p-3}\mod p&\t{if $p\e 1\mod 4$,}
\\-\f 1{16}\Gamma_p\ls 14^4\mod p&\t{if $p\e 3\mod 4$.}
\endcases$$
\pro{Theorem 2.1} Let $p$ be an odd prime.
\par $(\t{\rm i})$ If $p\e 3\mod 4$, then
$$A'_{\f{p-1}2}\e \f{p^2}3\b{\f{p-3}2}{\f{p-3}4}^{-2}\mod {p^3}.$$
\par $(\t{\rm ii})$ If $p\e 1\mod 4$ and so $p=x^2+4y^2$ with $x,y\in\Bbb Z$, then
$$A'_{\f{p-1}2}\e  4x^2-2p-\f{p^2}{4x^2}+3p^2x^2E_{p-3}+\f{p^2}2\sum_{k=1}^{(p-1)/2}\f{\b{2k}k^3}{64^k}
\Big(\sum_{i=1}^k\f 1{2i-1}\Big)^2\mod {p^3}.$$
\endpro
\par{\it Proof.} Since $A'_1=3$, the result is true for $p=3$. Now assume that $p>3$. For $\f p2<k<p$ we see that $p\mid \b{2k}k$. Thus, from Lemmas 2.5-2.7 we have
$$\align &\sum_{k=1}^{(p-1)/2}\f{\b{2k}k^3}{64^k}\sum_{i=1}^k\f 1{2i-1}
\\&\e\cases 0\mod {p^2}&\t{if $4\mid p-1$,}
\\-\f p{12}\cdot 16\cdot 2^{p-3}\b{(p-3)/2}{(p-3)/4}^{-2}\e -\f p3\b{(p-3)/2}{(p-3)/4}^{-2}\mod {p^2}&\t{if $4\mid p-3.$}
\endcases\endalign$$
and
$$\align&\sum_{k=1}^{(p-1)/2}\f{\b{2k}k^3}{64^k}\sum_{i=1}^k\f 1{(2i-1)^2}
\\&\e\cases \f 12(-4x^2)E_{p-3}=-2x^2E_{p-3}\mod p&\t{if $p\e 1\mod 4$,}
\\-2^{p-3}\b{(p-3)/2}{(p-3)/4}^{-2}\e -\f 14\b{(p-3)/2}{(p-3)/4}^{-2}\mod p&\t{if $p\e 3\mod 4$.}
\endcases\endalign$$
By Lemmas 2.1-2.2,
$$\sum_{k=0}^n\b nk\b{n+k}k\f{\b{2k}k}{(-4)^k}D_k=0\q\t{for $n=1,3,5,\ldots$}.\tag 2.1$$
Note that $\b nk\b{n+k}k=\b{2k}k\b{n+k}{2k}$. By [Su2, Lemma 2.2],
$$\b{\f{p-1}2+k}{2k}\e \f{\b{2k}k}{(-16)^k}\mod {p^2}\qtq{for}k=1,2,\ldots,\f{p-1}2.\tag 2.2$$
Hence, for $p\e 3\mod 4$,
$$\sum_{k=0}^{(p-1)/2}\f{\b{2k}k^3}{64^k}D_k\e \sum_{k=0}^{(p-1)/2}
\b {\f{p-1}2}k\b{\f{p-1}2+k}k\f{\b{2k}k}{(-4)^k}D_k=0\mod {p^2}.$$
\par Now, from the above and Lemmas 2.3-2.4
 we deduce that for $p\e 3\mod 4$,
$$\align A'_{\f{p-1}2}&\e1+\sum_{k=1}^{(p-1)/2}\f{\b{2k}k^3}{64^k}
\Big(1-p\sum_{i=1}^k\f 1{2i-1}+\f{p^2}2\Big(D_k-2\sum_{i=1}^k\f 1{(2i-1)^2}\Big)\Big)
\\&\e -\f{p^2}4\b{(p-3)/2}{(p-3)/4}^{-2}+\f{p^2}3\b{(p-3)/2}{(p-3)/4}^{-2}
+\f{p^2}4\b{(p-3)/2}{(p-3)/4}^{-2}
\\&=\f{p^2}3\b{(p-3)/2}{(p-3)/4}^{-2}\mod{p^3},\endalign$$
and for $p=x^2+4y^2\e 1\mod 4$,
$$\align A'_{\f{p-1}2}&\e1+\sum_{k=1}^{(p-1)/2}\f{\b{2k}k^3}{64^k}
\Big(1-p\sum_{i=1}^k\f 1{2i-1}
\\&\q+\f{p^2}2\Big(\Big(\sum_{i=1}^k\f 1{2i-1}\Big)^2-3\sum_{i=1}^k\f 1{(2i-1)^2}\Big)\Big)
\\&\e 4x^2-2p-\f{p^2}{4x^2}+\f{p^2}2\sum_{k=1}^{(p-1)/2}\f{\b{2k}k^3}{64^k}
\Big(\sum_{i=1}^k\f 1{2i-1}\Big)^2\\&\q-\f 32p^2(-2x^2E_{p-3})
\mod {p^3}.\endalign$$
This completes the proof.
\v2
\pro{Conjecture 2.1} Let $p$ be a prime of the form $4k+1$ and so $p=x^2+4y^2$
with $x,y\in\Bbb Z$. Then
$$\sum_{k=1}^{(p-1)/2}\f{\b{2k}k^3}{64^k}
\Big(\sum_{i=1}^k\f 1{2i-1}\Big)^2\e \f 23x^2E_{p-3}\mod p$$
and so
$$A'_{\f{p-1}2}\e 4x^2-2p-\f{p^2}{4x^2}+\f{10}3p^2x^2E_{p-3}\mod {p^3}.$$
\endpro

\v2
\par We point out the author's following conjectures in [Su7].
\v2
\par{\bf Conjecture 2.2} Let $p>3$ be a prime and $m,r\in\Bbb Z^+$. Then
$$\align &A'_{mp^r-1}-A'_{mp^{r-1}-1}\\&\e\f 53m^3\Big(\sum_{k=1}^m\b mk\b{m-1}{k-1}\b{m+k-1}{k-1}\Big)p^{3r}B_{p-3}\mod {p^{3r+1}}.\endalign$$
\par{\bf Conjecture 2.3} Let $p>3$ be a prime,  $m,r\in\Bbb Z^+$ and $C'_m$ be given in (1.2). Then
$$A'_{mp^r}\e A'_{mp^{r-1}}+C'_mp^{3r}\Big(\f{B_{2p-4}}{2p-4}-2\f{B_{p-3}}{p-3}\Big)\mod {p^{3r+2}}.$$

\par{\bf Conjecture 2.4} Let $p>5$ be a prime, $m,r\in\Bbb Z^+$ and $C_m$ be given in (1.1). Then
$$ A_{mp^r}-A_{mp^{r-1}}\e 2C_mp^{3r}\Big(\f{B_{2p-4}}{2p-4}-2\f{B_{p-3}}{p-3}\Big)\mod {p^{3r+2}}.$$

\par{\bf Conjecture 2.5} Let $p>3$ be a prime and $m,r\in\Bbb Z^+$. Then there exists an odd number $c_m$ depending only on $m$ such that
$$A_{mp^r-1}-A_{mp^{r-1}-1}\e \f 23m^3c_mp^{3r}B_{p-3}\mod {p^{3r+1}}.$$
Moreover, $c_1=1,\; c_2=1,\; c_3=-17,\; c_4=-703,\; c_5=-21499,\; c_6=-628145.$

\section*{3. Identities and congruences for $t_n$}

\par\q\ Let $\{t_n\}$ be defined by (1.5).
\pro{Theorem 3.1} For $n=0,1,2,\ldots$ we have
$$t_n=(2n+1)!\sum_{k=0}^n\f {\b{2k}k}{4^k(2(n-k)+1)}.$$
\endpro
\par{\it Proof.} For $|x|<1$ it is well known that
$$\align&\t{\rm arctanh}(x)=\f 12\ln\f{1+x}{1-x}=\sum_{m=0}^{\infty}\f{x^{2m+1}}{2m+1},
\\&\f 1{\sqrt{1-x^2}}=\sum_{k=0}^{\infty}\b{-\f 12}k(-x^2)^k=
\sum_{k=0}^{\infty}\f{\b{2k}k}{4^k}x^{2k}.\endalign$$
From [Sl, A028353] we know that
$$\f{\t{\rm arctanh}(x)}{\sqrt{1-x^2}}=\sum_{n=0}^{\infty}t_n\f{x^{2n+1}}{(2n+1)!}
\q(|x|<1).$$
Thus,
$$\sum_{n=0}^{\infty}t_n\f{x^{2n+1}}{(2n+1)!}
=\Big(\sum_{m=0}^{\infty}\f{x^{2m+1}}{2m+1}\Big)
\Big(\sum_{k=0}^{\infty}\f{\b{2k}k}{4^k}x^{2k}\Big).$$
Comparing the coefficients of $x^{2n+1}$ on both sides yields the result.
\v2
\pro{Theorem 3.2}  For $n=0,1,2,\ldots$ we have
$$\align t_n^2&=-(2n+1)!^2\sum_{k=0}^{2n+1}\b{2n+1+k}{2k}\b{2k}k^2\f 1{(-4)^k}
\sum_{i=1}^k\f 1{(2i-1)^2}
\\&=-(2n+1)!^2\sum_{k=0}^{2n+1}\b{2n+1+k}{2k}\b{2k}k^2\f 1{(-4)^k}
\Big(\sum_{i=1}^k\f 1{2i-1}\Big)^2.\endalign$$
\endpro
\par{\it Proof.} Let
$$\align &S_1(n)=\sum_{k=0}^{2n+1}\b{2n+1+k}{2k}\b{2k}k^2\f 1{(-4)^k}
\sum_{i=1}^k\f 1{(2i-1)^2},
\\&S_2(n)=-\Big(\sum_{k=0}^n\f {\b{2k}k}{4^k(2(n-k)+1)}\Big)^2
.\endalign$$
It is easy to see that
$$\align &S_1(0)=-1=S_2(0),\q S_1(1)=-\f{25}{36}=S_2(1),\\& S_1(2)=-\Ls {89}{120}^2=S_2(2),\q S_1(3)=-\Ls{381}{560}^2=S_2(3).\endalign$$
Using the Maple software doublesum.mpl
and the method in [CHM], we find that for $i=1,2$ and $n=0,1,2,\ldots$,
$$\align
&4(n+4)^2(2n+7)^2(2n+9)^2(4n+9)(75+72n+16n^2)S_i(n+4)
\\&-(2n+7)^2(6913575+17355348n+18370228n^2 +10658464n^3 +3670400 n^4
\\&\q+751872n^5 + 84992n^6 + 4096n^7)S_i(n+3)
\\&+(4n+11)(18889425+56173260n+72583012n^2+53324832n^3+24399376n^4
\\&\q+7128000n^5 +1299328n^6+135168 n^7+ 6144n^8)S_i(n+2)
\\&-8(n+2)^2(1254375+3543600n+4277038n^2+2861712n^3+1146240n^4
\\&\q+274560n^5+36352n^6 +2048n^7)S_i(n+1)
\\&+16(n+1)^2(n+2)^2(2n+3)^2(4n+13)(163+104n+16n^2)S_i(n) =0.
\endalign$$
Hence, for $n=0,1,2,\ldots$ we have $S_1(n)=S_2(n)$. That is,
$$\Big(\sum_{k=0}^n\f {\b{2k}k}{4^k(2(n-k)+1)}\Big)^2=
-\sum_{k=0}^{2n+1}\b{2n+1+k}{2k}\b{2k}k^2\f 1{(-4)^k}
\sum_{i=1}^k\f 1{(2i-1)^2}.$$
This together with Theorem 3.1 and (2.1) yields the result.
\v2
\pro{Theorem 3.3} Let $p$ be an odd prime. Then
$$\align &t_p\e \big(1+4(-1)^{\f{p-1}2}\big)p^2\mod {p^3},
\\&t_{p-1}\e (-1)^{\f{p-1}2}\big(2p+2^p-2+(pB_{p-1})^2\big)\mod {p^2},
\\& t_{\f{p-1}2}\e  pB_{p-1}-p+2^{p-1}-1\mod {p^2},
\\& t_{\f{p+1}2}\e  pB_{p-1}-3p+2^{p-1}-1\mod {p^2}
\endalign$$
and
$$t_{\f{p-3}4}\e \pm \f 1{\b{(p-1)/2}{(p-3)/4}}\mod p\q\t{for $p\e 3\mod 4$}.$$
\endpro
\par{\it Proof.} Since  $p\mid \b{2k}k$ for $\f p2<k<p$ and
$$\f{(2p+1)!}{p^2}=(p-1)!(p+1)\cdots(2p-1)\cdot 2(2p+1)\e 2(p-1)!^2\e 2\mod p,$$ we see that
$$\align t_p&=(2p+1)!\sum_{k=0}^p\f{\b{2k}k}{4^k(2(p-k)+1)}
\\&\e (2p+1)!\Big(\sum_{k=0}^{(p-1)/2}\f{\b{2k}k}{4^k(2(p-k)+1)}
+\f{\b{p+1}{(p+1)/2}}{4^{\f{p+1}2}\cdot p}+\f{\b{2p}p}{4^p}\Big)
\\&= (2p+1)!\Big(\sum_{k=0}^{(p-1)/2}\f{\b{2k}k}{4^k(2p+1-2k)}+
\f{\b{p-1}{(p-1)/2}}{4^{\f{p-1}2}(p+1)}+\f{\b{2p-1}{p-1}}{2\cdot 4^{p-1}}\Big)
\\&\e 2p^2\sum_{k=0}^{(p-1)/2}\f{\b{2k}k}{4^k(p+1-2k)}+2p^2\cdot(-1)^{\f{p-1}2}+p^2
\\&=p^2\sum_{k=0}^{(p-1)/2}\b{-1/2}k(-1)^k\f 1{\f {p+1}2-k}+2p^2(-1)^{\f{p-1}2}+p^2
\\&\e p^2\sum_{k=0}^{(p-1)/2}\b{(p-1)/2}k(-1)^k\f 1{\f {p+1}2-k}+2(-1)^{\f{p-1}2}p^2+p^2
\mod {p^3}.
\endalign$$
By [G, (1.43)],
$$\sum_{k=0}^n\b nk(-1)^k\f 1{x-k}=\f{(-1)^n}{(x-n)\b xn}.\tag 3.1$$
Hence
$$ \sum_{k=0}^{(p-1)/2}\b{(p-1)/2}k(-1)^k\f 1{\f {p+1}2-k}
= \f{(-1)^{\f{p-1}2}}{(\f{p+1}2-\f{p-1}2)\b{(p+1)/2}{(p-1)/2}}\e 2(-1)^{\f{p-1}2}\mod p.$$
Therefore,
$$t_p\e 2(-1)^{\f{p-1}2}p^2+2(-1)^{\f{p-1}2}p^2+p^2=\big(1+4(-1)^{\f{p-1}2}\big)
p^2\mod {p^3}.$$
\par Since $p\mid (2p-1)!$ and $p\mid \b{2k}k$ for $\f p2<k<p$  we see that
$$\align  t_{p-1}&=(2p-1)!\sum_{k=0}^{p-1}\f{\b{2k}k}{4^k(2p-1-2k)}
\\&\e (2p-1)!\sum_{k=0}^{(p-3)/2}\f{\b{2k}k}{4^k(2p-1-2k)}+(2p-1)!\f{\b{p-1}{(p-1)/2}}
{4^{\f{p-1}2}\cdot p}
\\&\e\f{(2p-1)!}2\sum_{k=0}^{(p-3)/2}\f{\b{2k}k}{(-4)^k}(-1)^k\f 1{\f{2p-1}2-k}+(p-1)!^2\f{\b{p-1}{(p-1)/2}}{2^{p-1}}
\mod {p^2}.\endalign$$
It is well known (see for example [Su1]) that $H_{\f{p-1}2}\e -2q_p(2)\mod p$ and $(p-1)!\e pB_{p-1}-p\mod {p^2}$. Thus,
$$\align \f{\b{p-1}{(p-1)/2}}{2^{p-1}}&=\f{(p-1)(p-2)\cdots(p-\f{p-1}2)}
{(\f{p-1}2)!\cdot 2^{p-1}}
\e (-1)^{\f{p-1}2}\f{1-pH_{\f{p-1}2}}{2^{p-1}}
\\&\e (-1)^{\f{p-1}2}\f{(1+pq_p(2))^2}{2^{p-1}}=(-1)^{\f{p-1}2}2^{p-1}\mod {p^2}\endalign$$
and so
$$\align (p-1)!^2\f{\b{p-1}{(p-1)/2}}{2^{p-1}}&\e (-1)^{\f{p-1}2}(p-1)!^2(2^{p-1}-1+1)\\&\e (-1)^{\f{p-1}2}(2^{p-1}-1)+(-1)^{\f{p-1}2}(p-1)!^2
\\&\e (-1)^{\f{p-1}2}\big(2^{p-1}-1+(pB_{p-1}-p)^2\big)
\\&\e (-1)^{\f{p-1}2}\big(2^{p-1}-1+(pB_{p-1})^2+2p\big)\mod {p^2}.
\endalign$$
Since
$\f{\b{2k}k}{(-4)^k}=\b{-\f 12}k\e \b{\f{p-1}2}k\mod {p}$ for $k\le \f{p-1}2$, using (3.1) we see that
$$\align &\f{(2p-1)!}2\sum_{k=0}^{(p-3)/2}\f{\b{2k}k}{(-4)^k}(-1)^k\f 1{\f{2p-1}2-k}
\\&\e\f{(2p-1)!}2\sum_{k=0}^{(p-3)/2}\b{\f{p-1}2}k(-1)^k\f 1{\f{2p-1}2-k}
\\&=\f{(2p-1)!}2\sum_{k=0}^{(p-1)/2}\b{\f{p-1}2}k(-1)^k\f 1{\f{2p-1}2-k}
\\&\q-(-1)^{\f{p-1}2}(p^2-1^2)\cdots(p^2-(p-1)^2)
\\&=\f{p\cdot (p^2-1^2)\cdots(p^2-(p-1)^2)}2\cdot\f{(-1)^{\f{p-1}2}}
{(\f{2p-1}2-\f{p-1}2)\b{(2p-1)/2}{(p-1)/2}}
\\&\q-(-1)^{\f{p-1}2}(p^2-1^2)\cdots(p^2-(p-1)^2)
\\&\e (-1)^{\f{p-1}2}(p-1)!^2\f 1{\f{(\f p2+1)(\f p2+2)\cdots(\f p2+\f{p-1}2)}{(\f{p-1}2)!}}- (-1)^{\f{p-1}2}(p-1)!^2
\\&\e (-1)^{\f{p-1}2}(p-1)!^2\f 1{1+\f p2H_{\f{p-1}2}}- (-1)^{\f{p-1}2}(p-1)!^2
\\&\e (-1)^{\f{p-1}2}(p-1)!^2\Big(\f 1{1-pq_p(2)}-1\Big)
\\&\e (-1)^{\f{p-1}2}(p-1)!^2pq_p(2)\e (-1)^{\f{p-1}2}(2^{p-1}-1)
\mod {p^2}.\endalign$$
Therefore,
$$t_{p-1}\e (-1)^{\f{p-1}2}(2(2^{p-1}-1)+(pB_{p-1})^2+2p\big)\mod {p^2}.$$
\par From [T1],
$$\sum_{k=1}^{(p-1)/2}\f{\b{2k}k}{4^kk}\e \sum_{k=1}^{p-1}\f{\b{2k}k}{4^kk}
\e 2q_p(2)\mod p.$$
Thus,
$$\align t_{\f{p-1}2}&=p!\sum_{k=0}^{(p-1)/2}\f{\b{2k}k}{4^k(p-2k)}
\e(p-1)!\Big(1-\f p2\sum_{k=1}^{(p-1)/2}\f{\b{2k}k}{4^kk}\Big)
\\&\e (p-1)!(1-pq_p(2))\e (pB_{p-1}-p)(1-pq_p(2))
\\&\e pB_{p-1}-p+2^{p-1}-1\mod {p^2}.\endalign$$
From (1.5) and the above congruence for $t_{\f{p-1}2}$ modulo $p^2$ we see that
$$\align t_{\f{p+1}2}&\e \Big(8\Ls{p-1}2^2+12\Ls{p-1}2+5\Big)t_{\f{p-1}2}
\e (2p+1)t_{\f{p-1}2}
\\&\e (2p+1)\big(pB_{p-1}-p+2^{p-1}-1\big)\e pB_{p-1}-3p+2^{p-1}-1\mod {p^2}.\endalign$$
\par For $p\e 3\mod 4$, from Theorem 3.2 and (2.2) we see that
$$\align t_{\f{p-3}4}^2&=-\ls{p-1}2!^2
\sum_{k=0}^{(p-1)/2}\b{\f{p-1}2+k}{2k}\b{2k}k^2\f 1{(-4)^k}\sum_{i=1}^k\f 1{(2i-1)^2}
\\&\e -\ls{p-1}2!^2
\sum_{k=0}^{(p-1)/2}\f{\b{2k}k^3}{64^k}\sum_{i=1}^k\f 1{(2i-1)^2}
\mod {p^2}.\endalign$$
Since $p\e 3\mod 4$ we have $\ls{p-1}2!^2\e -(p-1)!\e 1\mod p.$
By the proof of Theorem 2.1,
$$\sum_{k=0}^{(p-1)/2}\f{\b{2k}k^3}{64^k}\sum_{i=1}^k\f 1{(2i-1)^2}
\e -\f 1{4\b{(p-3)/2}{(p-3)/4}^2}\mod p.$$
Thus,
$$t_{\f{p-3}4}^2\e \f 1{4\b{(p-3)/2}{(p-3)/4}^2}\mod p
\q\t{and so}\q t_{\f{p-3}4}\e \mp \f 1{2\b{(p-3)/2}{(p-3)/4}}\e \pm \f 1{\b{(p-1)/2}{(p-3)/4}}\mod p.$$
This completes the proof.

\end{document}